\documentclass[reqno]{amsart}
\usepackage[T1]{fontenc}
\usepackage{geometry}
\usepackage{mathtools}
\usepackage{amstext}
\usepackage{amsthm}
\usepackage{amssymb}
\usepackage{stmaryrd}

\makeatletter
\numberwithin{equation}{section}
\numberwithin{figure}{section}
\theoremstyle{plain}
\newtheorem{thm}{\protect\theoremname}
\theoremstyle{definition}
\newtheorem{defn}[thm]{\protect\definitionname}

\usepackage{shuffle}

\makeatother

\providecommand{\definitionname}{Definition}
\providecommand{\theoremname}{Theorem}

\begin{document}
\address[Minoru Hirose]{Graduate School of Science and Engineering, Kagoshima University, 1-21-35 Korimoto, Kagoshima, Kagoshima 890-0065, Japan}
\email{mathhirose@gmail.com}
\subjclass[2010]{11M32}
\keywords{multiple zeta values, parity result, multitangent function}
\title[An explicit parity theorem for MZVs]{An explicit parity theorem for multiple zeta values via multitangent
functions}
\author{Minoru Hirose}
\begin{abstract}
We give an explicit formula for the well-known parity result for multiple
zeta values as an application of the multitangent functions.
\end{abstract}

\maketitle

\section{Introduction}

Multiple zeta values are real numbers defined by
\[
\zeta(k_{1},\dots,k_{d}):=\sum_{0<m_{1}<\cdots<m_{d}}\frac{1}{m_{1}^{k_{1}}\cdots m_{d}^{k_{d}}}
\]
for positive integers $k_{1},\dots,k_{d}$ with $k_{d}\geq2$. Here,
$k_{1}+\cdots+k_{d}$ is called the \emph{weight,} and $d$ is called
the \emph{depth}. The following theorem is conjectured in \cite{BG-conj}
and proved or generalized by several authors (\cite{TsumuraParity}
the first proof, \cite{Brown-depth,IKZ,Machide} via double shuffle
relations, \cite{Jarossy} for associators, \cite{Panzer} for multiple
polylogarithms).
\begin{thm}[Parity result for multiple zeta values]
\label{thm:parity}When the weight $k_{1}+\cdots+k_{d}$ and the
depth $d$ are of the opposite parity, $\zeta(k_{1},\dots,k_{d})$
is a $\mathbb{Q}[\pi^{2}]$-linear combination of multiple zeta values
of depth at most $d-1$.
\end{thm}

The purpose of this paper is to give a new proof of the theorem and
provide an explicit formula by using a theory of multitangent functions
by Bouillot \cite{Bou14AlgMTGF}. Note that Kusunoki-Nakamura-Sasaki
\cite{KNS} also gave an explicit formula for Theorem \ref{thm:parity}.

For positive integers $k_{1},\dots,k_{d}$, let $\zeta^{*}(k_{1},\dots,k_{d})$
be the stuffle regularized multiple zeta values (see \cite{IKZ}),
\[
\zeta^{\star,*}(k_{1},\dots,k_{d}):=\sum_{\substack{\boxempty_{j}=\,,\,\text{or }+\\
\text{for }j=1,\dots,d-1
}
}\zeta^{*}(k_{1}\boxempty_{1}\cdots\boxempty_{d-1}k_{d})
\]
(e.g. $\zeta^{\star,*}(5,3,1)=\zeta^{*}(5,3,1)+\zeta^{*}(5,4)+\zeta^{*}(8,1)+\zeta^{*}(9)$)
the stuffle regularized multiple zeta star values. For $k_{1},\dots,k_{d}\in\mathbb{Z}_{\geq1}$
and $a\in\mathbb{Z}_{\geq0}$, we put
\[
\zeta_{a}^{*}(k_{1},\dots,k_{d}):=(-1)^{a}\sum_{\substack{a_{1}+\cdots+a_{d}=a\\
a_{1},\dots,a_{d}\geq0
}
}\zeta_{a}^{*}(k_{1}+a_{1},\dots,k_{d}+a_{d})\prod_{j=1}^{d}{k_{j}-1+a_{j} \choose a_{j}}.
\]
When $k_{d}\geq2$, we simply write $\zeta^{\star}(k_{1},\dots,k_{d})$
and $\zeta_{a}(k_{1},\dots,k_{d})$ instead of $\zeta^{\star,*}(k_{1},\dots,k_{d})$
and $\zeta_{a}^{*}(k_{1},\dots,k_{d})$, respectively. The following
is the main theorem of this paper.
\begin{thm}
\label{thm:main}For positive integers $k_{1},\dots,k_{d}$ with $k_{d}\geq2$
and $\sum_{j=1}^{d}k_{j}\not\equiv d\pmod{2}$, we have
\begin{align*}
\zeta(k_{1},\dots,k_{d}) & =\frac{\zeta(k_{1},\dots,k_{d})-\zeta^{\star}(k_{1},\dots,k_{d})}{2}\\
 & \quad\quad-\sum_{0\leq i<j\leq d}\sum_{\substack{a+2m+b=k_{j}\\
a,b,m\geq0
}
}(-1)^{m+i+b+k_{1}+\cdots+k_{j}}\frac{(2\pi)^{2m}}{2(2m)!}B_{2m}\\
 & \qquad\qquad\qquad\times\zeta^{\star,*}(k_{1},\dots,k_{i})\zeta_{a}^{*}(k_{j-1},\dots,k_{i+1})\zeta_{b}(k_{j+1},\dots,k_{d}).
\end{align*}
\end{thm}

Theorem \ref{thm:main} implies Theorem \ref{thm:parity} since both
$\zeta(k_{1},\dots,k_{d})-\zeta^{\star}(k_{1},\dots,k_{d})$ and $\zeta^{\star,*}(k_{1},\dots,k_{i})\times\zeta_{a}^{*}(k_{j-1},\dots,k_{i+1})\times\zeta_{b}(k_{j+1},\dots,k_{d})$
can be written as a $\mathbb{Q}$-linear sum of multiple zeta values
of depth $\leq d-1$.

\section{A proof}

For the proof of the main theorem, we use the theory of multitangent
functions due to Bouillot \cite{Bou14AlgMTGF}. Let us start with
the definition of the multitangent functions. For positive integers
$k_{1},\dots,k_{d}$ with $k_{1},k_{d}\geq2$, the multitangent tangent
functions are holomorphic functions of $z\in\mathbb{C}\setminus\mathbb{Z}$
defined by the absolutely convergent series
\[
\Psi_{k_{1},\dots,k_{d}}(z):=\sum_{-\infty<m_{1}<\cdots<m_{d}<\infty}\frac{1}{(z+m_{1})^{k_{1}}\cdots(z+m_{d})^{k_{d}}}.
\]
By the notation in \cite{Bou14AlgMTGF}, $\Psi_{k_{1},\dots,k_{d}}(z)$
is written as $\mathcal{T}e^{k_{d},\dots,k_{1}}$. Then, by letting
\[
i:=\max\{0\leq i\leq d\,\mid\,m_{i}<0\},\ j:=\min\{0\leq j\leq d\,\mid\,m_{j+1}>0\},\qquad(m_{0}:=-\infty,m_{d+1}:=\infty),
\]
$\Psi_{k_{1},\dots,k_{d}}$ is decomposed as
\begin{align*}
\Psi_{k_{1},\dots,k_{d}}(z) & =\sum_{0\leq i\leq j\leq d}\left(\sum_{-\infty<m_{1}<\cdots<m_{i}<0}\frac{1}{(z+m_{1})^{k_{1}}\cdots(z+m_{i})^{k_{i}}}\right)\\
 & \qquad\times\left(\sum_{0<m_{j+1}<\cdots<m_{d}<\infty}\frac{1}{(z+m_{j+1})^{k_{j+1}}\cdots(z+m_{d})^{k_{d}}}\right)\times\begin{cases}
1 & j-i=0\\
\frac{1}{z^{k_{j}}} & j-i=1\\
0 & j-i\geq2
\end{cases}\\
 & =\sum_{0\leq i\leq j\leq d}(-1)^{k_{1}+\cdots+k_{i}}\zeta^{(-z)}(k_{i},\dots,k_{1})\zeta^{(z)}(k_{j+1},\dots,k_{d})\times\begin{cases}
1 & j-i=0\\
\frac{1}{z^{k_{j}}} & j-i=1\\
0 & j-i\geq2
\end{cases},
\end{align*}
where
\[
\zeta^{(z)}(l_{1},\dots,l_{r})\coloneqq\sum_{0<n_{1}<\cdots<n_{r}}\frac{1}{(z+n_{1})^{l_{1}}\cdots(z+n_{r})^{l_{r}}}\qquad(l_{r}>1)
\]
is the Hurwitz multiple zeta values. For positive integers $l_{1},\dots,l_{r}$
with $l_{r}\geq2$, the Taylor expansion of the Hurwitz multiple zeta
value near $z=0$ is given by
\[
\zeta^{(z)}(l_{1},\dots,l_{r})=\sum_{a=0}^{\infty}z^{a}\zeta_{a}(l_{1},\dots,l_{r}).
\]
Based on this observation, define the stuffle regularized Hurwitz
multiple zeta values for positive integers $l_{1},\dots,l_{r}$ (not
necessarily $l_{r}\geq2$) by
\begin{equation}
\zeta^{(z),*}(l_{1},\dots,l_{r}):=\sum_{a=0}^{\infty}z^{a}\zeta_{a}^{*}(l_{1},\dots,l_{r}).\label{eq:streg_hurwitz}
\end{equation}
Then the series (\ref{eq:streg_hurwitz}) absolutely converges for
$\left|z\right|<1$ and is analytically continued to $z\in\mathbb{C}\setminus\mathbb{Z}$
(see \cite[Section 7]{Bou14AlgMTGF}).
\begin{defn}[{\cite[Section 7]{Bou14AlgMTGF}}]
For positive integers $k_{1},\dots,k_{d}$, the stuffle regularized
multitangent function is defined by
\[
\Psi_{k_{1},\dots,k_{d}}(z):=\sum_{0\leq i\leq j\leq d}(-1)^{k_{1}+\cdots+k_{i}}\zeta^{(-z),*}(k_{i},\dots,k_{1})\zeta^{(z),*}(k_{j+1},\dots,k_{d})\times\begin{cases}
1 & j-i=0\\
\frac{1}{z^{k_{j}}} & j-i=1\\
0 & j-i\geq2
\end{cases},
\]
or more explicitly,
\begin{align*}
\Psi_{k_{1},\dots,k_{d}}(z) & :=\sum_{j=0}^{d}\sum_{a=0}^{\infty}\sum_{b=0}^{\infty}(-1)^{k_{1}+\cdots+k_{j}+a}z^{a+b}\zeta_{a}^{*}(k_{j},\dots,k_{1})\zeta_{b}^{*}(k_{j+1},\dots,k_{d})\\
 & \quad+\sum_{j=1}^{d}\sum_{a=0}^{\infty}\sum_{b=0}^{\infty}(-1)^{k_{1}+\cdots+k_{j-1}+a}z^{a+b-k_{j}}\zeta_{a}^{*}(k_{j-1},\dots,k_{1})\zeta_{b}^{*}(k_{j+1},\dots,k_{d}).
\end{align*}
\end{defn}

Bouillot proved the following identity.
\begin{thm}[{Bouillot \cite[Theorem 6]{Bou14AlgMTGF}}]
\label{thm:Bou_mono_red}For positive integers $k_{1},\dots,k_{d}$,
we have
\[
\Psi_{k_{1},\dots,k_{d}}(z)=\delta^{k_{1},\dots,k_{d}}+\sum_{j=1}^{d}\sum_{\substack{a+s+b=k_{j}\\
a,b\geq0,\,s\geq1
}
}(-1)^{k_{1}+\cdots+k_{j-1}+a}\zeta_{a}^{*}(k_{j-1},\dots,k_{1})\zeta_{b}^{*}(k_{j+1},\dots,k_{d})\Psi_{s}^{*}(z)
\]
where
\[
\delta^{k_{1},\dots,k_{d}}:=\begin{cases}
\frac{(-1)^{n}\pi^{2n}}{(2n)!} & (k_{1},\dots,k_{d})=(\overbrace{1,\dots,1}^{2n})\ \text{for some }n\in\mathbb{Z}_{\geq0}\\
0 & \text{otherwise}.
\end{cases}
\]
\end{thm}

The coefficient of $z^{0}$ of the left-hand side of Theorem \ref{thm:Bou_mono_red}
is given by
\begin{align*}
 & \sum_{j=0}^{d}(-1)^{k_{1}+\cdots+k_{j}}\zeta^{*}(k_{j},\dots,k_{1})\zeta^{*}(k_{j+1},\dots,k_{d})\\
 & \quad+\sum_{j=1}^{d}\sum_{\substack{a,b\geq0\\
a+b=k_{j}
}
}(-1)^{k_{1}+\cdots+k_{j-1}+a}\zeta_{a}^{*}(k_{j-1},\dots,k_{1})\zeta_{b}^{*}(k_{j+1},\dots,k_{d}),
\end{align*}
and the coefficient of $z^{0}$ of the right-hand side of Theorem
\ref{thm:Bou_mono_red} is given by
\[
\delta^{k_{1},\dots,k_{d}}+\sum_{j=1}^{d}\sum_{\substack{a+2m+b=k_{j}\\
a,b\geq0,\,m\geq1
}
}(-1)^{k_{1}+\cdots+k_{j-1}+a}\zeta_{a}^{*}(k_{j-1},\dots,k_{1})\zeta_{b}^{*}(k_{j+1},\dots,k_{d})2\zeta(2m).
\]
Thus, we have
\begin{align}
 & \sum_{j=0}^{d}(-1)^{k_{1}+\cdots+k_{j}}\zeta^{*}(k_{j},\dots,k_{1})\zeta^{*}(k_{j+1},\dots,k_{d})\nonumber \\
 & \quad+\sum_{j=1}^{d}\sum_{\substack{a,b\geq0\\
a+b=k_{j}
}
}(-1)^{k_{1}+\cdots+k_{j-1}+a}\zeta_{a}^{*}(k_{j-1},\dots,k_{1})\zeta_{b}^{*}(k_{j+1},\dots,k_{d})\nonumber \\
 & =\delta^{k_{1},\dots,k_{d}}+\sum_{j=1}^{d}\sum_{\substack{a+2m+b=k_{j}\\
a,b\geq0,\,m\geq1
}
}(-1)^{k_{1}+\cdots+k_{j-1}+a}\zeta_{a}^{*}(k_{j-1},\dots,k_{1})\zeta_{b}^{*}(k_{j+1},\dots,k_{d})2\zeta(2m).\label{eq:fund_eq}
\end{align}
Since $B_{0}=1$ and
\[
\zeta(2m)=\frac{(-1)^{m+1}(2\pi)^{2m}}{2(2m)!}B_{2m},
\]
by multiplying $(-1)^{k_{1}+\cdots+k_{d}}$, the identity (\ref{eq:fund_eq})
can be rewritten as
\begin{align}
 & \sum_{j=0}^{d}(-1)^{k_{j+1}+\cdots+k_{d}}\zeta^{*}(k_{j},\dots,k_{1})\zeta^{*}(k_{j+1},\dots,k_{d})\nonumber \\
 & =\delta^{k_{1},\dots,k_{d}}+\sum_{j=1}^{d}\sum_{\substack{a+2m+b=k_{j}\\
a,b,m\geq0
}
}(-1)^{b+k_{j+1}+\cdots+k_{d}+m+1}\nonumber \\
 & \qquad\qquad\qquad\times\frac{(2\pi)^{2m}}{(2m)!}B_{2m}\zeta_{a}^{*}(k_{j-1},\dots,k_{1})\zeta_{b}^{*}(k_{j+1},\dots,k_{d}).\label{eq:fund_eq2}
\end{align}
Let $k_{1},\dots,k_{d}$ be positive integers. Let
\[
A(k_{1},\dots,k_{d})=\sum_{0\leq i\leq j\leq d}(-1)^{i+k_{j+1}+\cdots+k_{d}}\zeta^{\star,*}(k_{1},\dots,k_{i})\zeta^{*}(k_{j},\dots,k_{i+1})\zeta^{*}(k_{j+1},\dots,k_{d}).
\]
We will evaluate $A(k_{1},\dots,k_{d})$ in two different ways. Firstly,
we have
\begin{align}
 & A(k_{1},\dots,k_{d})\nonumber \\
 & =(-1)^{d}\zeta^{\star,*}(k_{1},\dots,k_{d})\nonumber \\
 & \quad+\sum_{i=0}^{d-1}(-1)^{i}\zeta^{\star,*}(k_{1},\dots,k_{i})\sum_{j=i}^{d}(-1)^{k_{j+1}+\cdots+k_{d}}\zeta^{*}(k_{j},\dots,k_{i+1})\zeta^{*}(k_{j+1},\dots,k_{d})\nonumber \\
 & =(-1)^{d}\zeta^{\star,*}(k_{1},\dots,k_{d})+\sum_{i=0}^{d-1}(-1)^{i}\zeta^{\star,*}(k_{1},\dots,k_{i})\delta^{k_{i+1},\dots,k_{d}}\nonumber \\
 & \quad+\sum_{i=0}^{d-1}(-1)^{i}\zeta^{\star,*}(k_{1},\dots,k_{i})\sum_{j=i+1}^{d}\sum_{\substack{a+2m+b=k_{j}\\
a,b,m\geq0
}
}(-1)^{b+k_{j+1}+\cdots+k_{d}+m+1}\nonumber \\
 & \qquad\qquad\times\frac{(2\pi)^{2m}}{(2m)!}B_{2m}\zeta_{a}^{*}(k_{j-1},\dots,k_{i+1})\zeta_{b}^{*}(k_{j+1},\dots,k_{d}).\label{eq:A1}
\end{align}
Here, we used (\ref{eq:fund_eq2}) for the second equality. Secondly,
we have
\begin{align}
A(k_{1},\dots,k_{d}) & =\sum_{j=0}^{d}(-1)^{k_{j+1}+\cdots+k_{d}}\zeta^{*}(k_{j+1},\dots,k_{d})\sum_{i=0}^{j}(-1)^{i}\zeta^{\star,*}(k_{1},\dots,k_{i})\zeta^{*}(k_{j},\dots,k_{i+1})\nonumber \\
 & =\sum_{j=0}^{d}(-1)^{k_{j+1}+\cdots+k_{d}}\zeta^{*}(k_{j+1},\dots,k_{d})\times\begin{cases}
1 & j=0\\
0 & j>0
\end{cases}\nonumber \\
 & =(-1)^{k_{1}+\cdots+k_{d}}\zeta^{*}(k_{1},\dots,k_{d}).\label{eq:A2}
\end{align}
Here, we used the well-known ``antipode identity'' 
\[
\sum_{i=0}^{j}(-1)^{i}\zeta^{\star,*}(k_{1},\dots,k_{i})\zeta^{*}(k_{j},\dots,k_{i+1})=0\qquad(j>0)
\]
for the second equality. By comparing (\ref{eq:A1}) and (\ref{eq:A2}),
we have the following theorem.
\begin{thm}
\label{thm:main2}For positive integers $k_{1},\dots,k_{d}\geq1$,
we have
\begin{align*}
 & (-1)^{d}\zeta^{\star,*}(k_{1},\dots,k_{d})-(-1)^{k_{1}+\cdots+k_{d}}\zeta^{*}(k_{1},\dots,k_{d})\\
 & =-\sum_{i=0}^{d-1}(-1)^{i}\zeta^{\star,*}(k_{1},\dots,k_{i})\delta^{k_{i+1},\dots,k_{d}}\\
 & \quad+\sum_{0\leq i<j\leq d}\sum_{\substack{a+2m+b=k_{j}\\
a,b,m\geq0
}
}(-1)^{i+b+k_{j+1}+\cdots+k_{d}+m}\frac{(2\pi)^{2m}}{(2m)!}B_{2m}\\
 & \qquad\times\zeta^{\star,*}(k_{1},\dots,k_{i})\zeta_{a}^{*}(k_{j-1},\dots,k_{i+1})\zeta_{b}^{*}(k_{j+1},\dots,k_{d}).
\end{align*}
\end{thm}

Theorem \ref{thm:main2} implies the following.

\begin{thm}
\label{thm:main3}For positive integers $k_{1},\dots,k_{d}\geq1$
with $\sum_{j=1}^{d}k_{j}\not\equiv d\pmod{2}$, we have
\begin{align*}
\zeta^{*}(k_{1},\dots,k_{d}) & =\frac{\zeta^{*}(k_{1},\dots,k_{d})-\zeta^{\star,*}(k_{1},\dots,k_{d})}{2}-\frac{1}{2}\sum_{i=0}^{d-1}(-1)^{d-i}\zeta^{\star,*}(k_{1},\dots,k_{i})\delta^{k_{i+1},\dots,k_{d}}\\
 & \quad-\sum_{0\leq i<j\leq d}\sum_{\substack{a+2m+b=k_{j}\\
a,b,m\geq0
}
}(-1)^{m+i+b+k_{1}+\cdots+k_{j}}\frac{(2\pi)^{2m}}{2(2m)!}B_{2m}\\
 & \qquad\qquad\times\zeta^{\star,*}(k_{1},\dots,k_{i})\zeta_{a}^{*}(k_{j-1},\dots,k_{i+1})\zeta_{b}^{*}(k_{j+1},\dots,k_{d}).
\end{align*}
\end{thm}

Theorem \ref{thm:main} is a special case of Theorem \ref{thm:main3}
since $\sum_{i=0}^{d-1}(-1)^{d-i}\zeta^{\star,*}(k_{1},\dots,k_{i})\delta^{k_{i+1},\dots,k_{d}}$
vanishes when $k_{d}\geq2$.

\subsection*{Acknowledgements}

The author thanks Ryota Umezawa and Shingo Saito for their useful
comments. This work was supported by JSPS KAKENHI Grant Numbers JP18K13392
and JP22K03244.

\end{document}